\documentclass[10pt]{amsart}
\usepackage{amsmath}
\usepackage{amssymb}
\usepackage{graphicx}
\usepackage{mathrsfs}
\usepackage{enumitem}
\usepackage{hyperref}
\usepackage{ytableau}
\usepackage{tikz}
\everymath{\displaystyle}
\usepackage{cite}

\newtheorem{theorem}{Theorem}
\newtheorem{lemma}[theorem]{Lemma}

\newtheorem{proposition}{Proposition} 
\newtheorem{definition}{Definition}[section]

\title[]{Expected value of statistics on type-B permutation tableaux}

\author[Althoff]{Ryan Althoff}
\author[Diethrich]{Daniel Diethrich}
\author[Lohss]{Amanda Lohss}
\author[Low]{Xin-Dee Low}
\author[Wichert]{Emily Wichert$\dagger$}

\address{Department of Computing, Mathematics and Statistics, Messiah University, One University Avenue, Mechanicsburg, 
	PA  17055, USA} 

\email{ra1241@messiah.edu}
\email{dd1278@messiah.edu}
\email{alohss@messiah.edu}
\email{lxindee@gmail.com}
\email{ewichert3@gmail.com}

\thanks{$\dagger$ This research was supported by NSF grant \#DMS-1722563}

\begin{document}

\begin{abstract}
Type-B permutation tableaux are combinatorial objects introduced by Lam and Williams that have an interesting connection with the partially asymmetric simple exclusion process (PASEP). In this paper, we compute the expected value of several statistics on these tableaux. Some of these computations are motivated by a similar paper on permutation tableaux. Others are motivated by the PASEP. In particular, we compute the expected number of rows, unrestricted rows, diagonal ones, adjacent south steps, and adjacent west steps. 
 \end{abstract}  

 
\maketitle

\section{Introduction}

Permutation tableaux are combinatorial objects introduced in~\cite{SW} that are in bijection with permutations~(see \cite{CN}[Theorem~11]).  These tableaux are also connected to the partially asymmetric simple exclusion process~\cite{CW2}, an important model in statistical mechanics. One particular variation of these tableaux are type-B permutation tableaux, introduced in~\cite{LW}, which are in bijection with signed permutations (see \cite{CJVW}[Theorem~4]). 

There have been several recent papers on permutation tableaux (see \cite{Bur}, \cite{CH}, \cite{CN}, \cite{CW3}, \cite{CW2}, and \cite{SW}). In particular, Corteel and Hitczenko calculated the expected value of several statistics on permutation tableaux in \cite{CH}. In this paper, we consider these same statistics on the type-B variation as well as other statistics that are interesting due to their connection with the PASEP. We will compute the expected value of these statistics using the approach developed in \cite{HL} for type-B permutation tableaux which is analogous to the approach used in \cite{CH} for permutation tableaux.

The paper is organized as follows. In Section~\ref{Prelim}, we introduce some preliminary definitions and terminology. In Section ~\ref{PASEP}, we discuss the relationship between permutation tableaux and the partially asymmetric simple exclusion process (PASEP) which is the motivation for several of our results. In Section ~\ref{PApproach}, we will discuss the probabilistic approach developed in \cite{HL} that we will use for our calculations.  The remaining sections include our results on the expected value of statistics on type-B permutation tableaux.

\section{Preliminaries \label{Prelim}} 
A \textit{Ferrers diagram} is a left aligned sequence of cells with weakly decreasing rows.  A Ferrers diagram is said to have half-perimeter \textit{$n$}, where \textit{$n$} equals the number of rows plus the number of columns. Each Ferrers diagram has a unique southeast border with the number of border edges equal to the half-perimeter. We will refer to these edges as \textit{steps} (south or west) with the first step beginning in the northeast corner of the diagram. We let $S_k$ denote a south step at the $k$th position and $W_k$ denote a west step.

\begin{figure}[h]
	\begin{center}
		\begin{ytableau}
		\none\\
		\none\\
		\none&\none&\none[
		\begin{picture} (3.5,7)
		\put(9.45,-5.0){\line(1,0){15}}
		\end{picture}
		]\\
		\none &  &  \\
		\none[
			\begin{picture} (3.5,7)
			\put(9.45,-4.6){\line(0,1){15}}
			\end{picture}
			] \\
		\none & \none & \none[(i)]
		\end{ytableau}
	\hspace{2cm}
		\begin{ytableau}
		\none & \\
		\none & & \\
		\none & & & \\
		\none &  &  \\
		\none[
		\begin{picture} (3.5,7)
		\put(9.45,-4.6){\line(0,1){15}}
		\end{picture}
		] \\
		\none & \none & \none[(ii)]
	\end{ytableau}
	\caption{(i) A Ferrers diagram with half-perimeter $5$ and steps $W_1S_2W_3W_4S_5$. (ii) A shifted Ferrers diagram obtained from (i) with the same half-perimeter and border edges.}
	\end{center}
\end{figure}
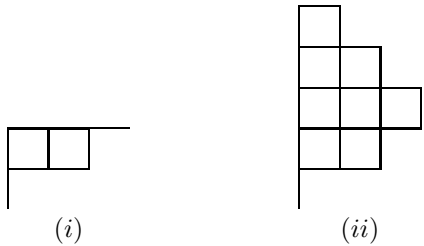

A \textit{shifted Ferrers diagram} is obtained from a Ferrers diagram by inserting additional rows above. If the Ferrers diagram has $k$ columns then the shifted Ferrers diagram is the same diagram but with $k$ rows inserted above of lengths $k$, $k-1$, $\dots$, 1, respectively. The rightmost cells of each added row are called \textit{diagonal cells}. A shifted Ferrers diagram has the same half-perimeter and border edges as that of the original diagram.

Permutation tableaux were introduced in \cite{SW} as fillings of Ferrers diagrams and type-B permutation tableaux, introduced in \cite{LW}, can be defined as fillings of shifted Ferrers diagrams (see \cite{CK}).

\begin{definition}
	\normalfont A \textit{permutation tableau} of size $n$ is a Ferrers diagram with half-perimeter $n$ in which each cell contains either a $0$ or a $1$ and satisfies the following conditions:
	\begin{enumerate}
		\item Each column must contain at least one $1$.
		\item A $0$ cannot have both a $1$ above it in the same column and a $1$ to left of it in the same row.
	\end{enumerate}
\end{definition}

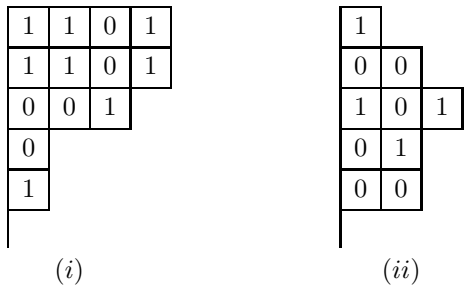
\begin{figure}[h]
	\begin{center}
		\begin{ytableau}
			\none & 1 & 1 & 0 & 1 \\
			\none & 1 & 1 & 0 & 1 \\
			\none & 0 & 0 & 1 \\
			\none & 0 \\
			\none & 1 \\
			\none[
			\begin{picture} (3.5,7)
			\put(9.45,-3.9){\line(0,1){14}}
			\end{picture}
			] \\
			\none & \none & \none[(i)]
		\end{ytableau}
		\hspace{2cm}
		\begin{ytableau}
			1 \\
			0 & 0 \\
			1 & 0 & 1 \\
			0 & 1 \\
			0 & 0 \\
			\none[
			\begin{picture} (3.5,7)
			\put(-5.95,-3.9){\line(0,1){14}}
			\end{picture}
			]\\
			\none & \none[(ii)]
		\end{ytableau}
		\caption{(i) A permutation tableau of size $n=10$. (ii) A type-B permutation tableau of size $n=6$.} \label{permEgs}
	\end{center}
\end{figure}

\begin{definition}\label{bdef}
	\normalfont A \textit{type-B permutation tableau} of size $n$ is a shifted Ferrers diagram with half-perimeter $n$, in which each cell contains either a $0$ or a $1$ and satisfies the following conditions:
	\begin{enumerate}
		\item Each column must contain at least one $1$.
		\item A $0$ cannot have both a $1$ above it in the same column and a $1$ to left of it in the same row.
		\item If a diagonal cell contains a $0$, then all of cells in that row must also be $0$.
	\end{enumerate}
\end{definition}

An important statistic for our calculations later is that of an \textit{unrestricted row}. A \textit{restricted $0$} is a $0$ that has a $1$ above it in the same column or a $0$ in a diagonal cell. Such a $0$ would require all cells to the left in the same row to be $0$'s since, by definition, a $0$ cannot have both a $1$ above it in the same column and a $1$ to the left of it in the same row. An unrestricted row is a row that contains no restricted $0$'s. 

We denote the number of unrestricted rows on the $k$th step by $U_k$. For example, consider the type-B permutation tableau in Figure~\ref{permEgs}(ii). This tableau has three unrestricted rows on the last step ($U_6=3$). Prior to that last step, the number of unrestricted rows are $U_1=1$, $U_2=2$, $U_3=3$, $U_4=2$, and $U_5=2$.


In this paper, we will calculate the following statistics on type-B permutation tableaux of size $n$,
\begin{enumerate}
	\item The expected number of rows (Theorem~\ref{thm:exprows}).
	\item The expected number of unrestricted rows (Theorem~\ref{thm:expunrows}).
	\item The expected number of $1$'s on the diagonal (Theorem~\ref{thm:exponesdiagonal}).
\end{enumerate}
These calculations are motivated by analogous work for permutation tableaux in \cite{CH}. In that paper, all three of those statistics were calculated for permutation tableaux as well as a fourth statistic which we were unable to calculate here for type-B permutation tableaux (see \cite{CH}[Theorem 1]).

We will also calculate the following statistics on type-B permutation tableaux of size $n$,
\begin{enumerate}
	\item The expected number of two adjacent south steps (Theorem~\ref{thm:twosouthstep}).
	\item The expected number of two adjacent west steps (Theorem~\ref{thm:twoadjacentw}).
\end{enumerate}
These statistics are interesting due to the connection between permutation tableaux and the PASEP. This will be discussed in detail in the next section, but these two statistics, along with a result from \cite{HL}[Theorem 6], provide the expected value of all possible configurations of two adjacent steps.

\section{The PASEP}\label{PASEP Section}
An important motivation for studying permutation tableaux is due to the connection between certain types of tableaux with the partially asymmetric simple exclusion process (PASEP). The PASEP is an important particle model in statistical mechanics that can be described as Markov chain, a stochastic process that transitions from one state to another with probabilities dependent only on the previous state. 

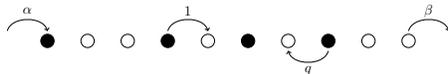
\begin{figure}[htbp] 
	\begin{center}
		\resizebox{6cm}{!}{
			\begin{tikzpicture}
			
			\node[draw,circle,fill=black] (c1) at (1,0) {};
			\node[draw,circle] (c2) at (2,0) {};
			\node[draw,circle] (c3) at (3,0) {};
			\node[draw,circle,fill=black] (c4) at (4,0) {};
			\node[draw,circle] (c5) at (5,0) {};
			\node[draw,circle,fill=black] (c6) at (6,0) {};
			\node[draw,circle] (c7) at (7,0) {};
			\node[draw,circle,fill=black] (c8) at (8,0) {};
			\node[draw,circle] (c9) at (9,0) {};
			\node[draw,circle] (c10) at (10,0) {};
			
			\node[-,black](1) at (.5,.75) {$\alpha$};
			\node[-,black](1) at (4.5,.75) {$1$};
			\node[-,black](1) at (7.5,-.75) {$q$};
			\node[-,black](1) at (10.5,.75) {$\beta$};

			\draw[->](4,.25) to [bend left=75] (5,.25) ;
			\draw[->](8,-.25) to [bend left=75] (7,-.25) ;
			\draw[->](10,.25) to [bend left=75] (11,.25);

			\draw[->](0,.25) to [bend left=75] (1,.25);

			\end{tikzpicture}	  
		}
	\end{center}
	\caption{An example of the PASEP as defined by a Markov chain of size $10$.}\label{PASEP}
\end{figure}

In particular, the PASEP is a Markov Chain on $n$ sites where each site is either empty or occupied by one particle. Empty sites are denoted by a $\circ$ and occupied sites by a $\bullet$. Particles enter and exit the system in one direction through the first and last sites. A particle enters from the left with rate $\alpha$ and exits the system from the right with rate $\beta$. Within the system, particles hop one site at a time with rate $1$ for a right hop and rate $q$ for a left hop. See Figure ~\ref{PASEP} below for an example of particular state of the PASEP when $n=10$.

Various types of tableaux have been used to give a combinatorial formula for the steady state distribution of the PASEP including permutation tableaux but also tree--like tableaux \cite{ABN}, alternative tableaux \cite{XV},and stammering tableaux \cite{J-V}. Each particular tableau is associated with a state of the PASEP which is called the tableau's {\it type}.

\begin{figure}[h]
 	\begin{center}
 		\begin{ytableau}
 			\none & 1 & 1 & 0 & 1 \\
 			\none & 0 & 0 & 1 \\
 			\none & 1 \\
 			\none[
 			\begin{picture} (3.5,7)
 			\put(9.45,-3.9){\line(0,1){14}}
 			\end{picture}
 			]\\
 			\none[
 			\begin{picture}(0,140)
 			
 			\put(14, 0){$\circ$}
 			\put(21, 0){$\bullet$}
 			\put(28, 0){$\circ$}
 			\put(35, 0){$\circ$}
 			\put(42, 0){$\bullet$}
 			\put(49, 0){$\circ$}
 			\put(56, 0){$\bullet$}
 			
 			\end{picture}

 			] \\
 			\none & \none & \none[(i)]
 		\end{ytableau}
 		\hspace{2cm}
 		\begin{ytableau}
 			1 \\
 			0 & 0 \\
 			0 & 1 \\
 			\none[
 			\begin{picture} (3.5,7)
 			\put(-5.95,-3.9){\line(0,1){14}}
 			\end{picture}
 			]\\
 			\none[
 			\begin{picture}(0,140)
 			
 			\put(-14, 0){$\bullet$}
 			\put(-7, 0){$\circ$}
 			\put(0, 0){$\circ$}
 			\put(7, 0){$\bullet$}
 			\put(14, 0){$\bullet$}
 			\put(21, 0){$\circ$}
 			\put(28, 0){$\circ$}
 			\put(35, 0){$\bullet$}
 			\end{picture}
 			] \\
 			\none & \none[(ii)]
 		\end{ytableau}
 		\caption{(i) A permutation tableau of size $8$ and its type, a Markov chain on $7$ sites. (ii) A type-B permutation tableau of size $4$ and its type, a symmetric Markov chain on $8$ sites.} \label{Tabtype}
 	\end{center}
 \end{figure}
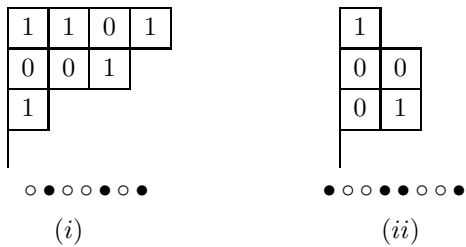

 The type of each permutation and type-B permutation tableau corresponds to the shape of its southeast border. In particular, each south step corresponds to a filled site and each west step corresponds to an empty site. For permutation tableaux, the first step is ignored since this step is always south. For type-B permutation tableaux, the Markov chain is doubled, preserving the site structure. Therefore, type-B permutation tableaux are connected with symmetric states of the PASEP. See Figure~\ref{Tabtype} for an example of the type associated with a particular permutation tableau and a type-B permutation tableau.

 In~\cite{HL}, the expected number of \textit{corners} in permutation and type-B permutation tableaux were computed (see Theorems 5 and 10). The motivation for these calculations was because a corner (south step followed by a west step) corresponds to a position in the PASEP where a filled site is followed by an empty site. From the number of corners, it is immediate to determine the number of \textit{inner corners} which is a west step followed by a south step.
 
 In this paper, we extend these results to provide calculations for all possible configurations of adjacent sites in the case of type-B permutation tableaux. We calculate the number of adjacent south steps in Section ~\ref{sectiontworows} which correspond to adjacent filled sites in the PASEP, and in Section~\ref{sectiontwocols}, we calculate the number of adjacent west steps which correspond to adjacent empty sites. This connection with the PASEP is the motivation for these two calculations  as mentioned in the introduction and is also further motivation for calculating the number of rows which corresponds to the number of filled sites in the PASEP.

\section{Probabilistic approach}\label{PApproach}
In \cite{HL}, probabilistic techniques were developed in order to compute the number of corners in type-B permutation tableaux. This approach is similar to the one developed in \cite{CH} for permutation tableaux. We will describe this approach here and use it to obtain our results in the subsequent sections. 

Let $\mathcal{B}_n$ denote the set of all type-B permutation tableaux of size $n$. Note that $|\mathcal{B}_n|=2^n\cdot n!$ as these tableaux are in bijection with signed permutations. Define $\mathbb{P}_n$ to be the probability distribution defined on $\mathcal{B}_n$ such that for all $T\in\mathcal{B}_n$,
\[
\mathbb{P}_n(T)=\frac{1}{2^n\cdot n!},
\]
and let $\mathbb{E}_n$ denote the expected value with respect to this measure.

There is a relationship between the measures $\mathbb{P}_{n-1}$ and $\mathbb{P}_n$ through an extension procedure. Any tableau of size $n-1$ can be extended to a tableau of size $n$ by inserting a new row or a new column. If a new column is inserted, there are several possible ways to fill the new column and these possibilities depend on the number of unrestricted rows. 

Recall that an unrestricted row does not contain a $0$ with a $1$ above it and therefore, this new column could place a $0$ or a $1$ in that row without contradicting the definition of type-B permutation tableaux. So as we extend tableaux of size $n-1$, the unrestricted rows are where choices can be made to determine the particular tableau of the extension. As mentioned in the preliminaries, we denote this key statistic by $U_k$ where $U_k=U_k(T)$ is the number of unrestricted rows on the $k$th step of a fixed tableau $T\in\mathbb{B}_n$.

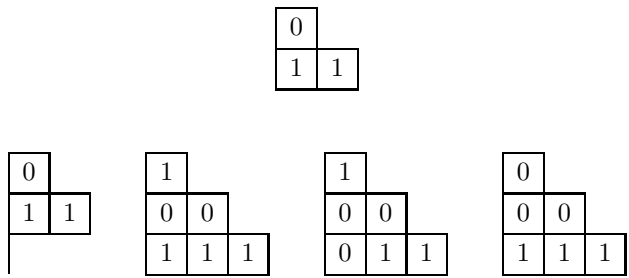
\begin{figure}[h] \label{extensions}
	\begin{center}
		\begin{ytableau}
			\none\\
			0 \\
			1&1 
		\end{ytableau}

\vspace{.25cm}
			\begin{ytableau}
			\none\\
			0 \\
			1&1 \\
			\none[
			\begin{picture} (3.5,7)
			\put(-6,-4.7){\line(0,1){15}}
			\end{picture}]&\none
		\end{ytableau}
	\hspace{.5cm}
		\begin{ytableau}
		\none\\
		1 \\
		0&0 \\
		1&1&1
	\end{ytableau}
	\hspace{.5cm}
\begin{ytableau}
	\none\\
	1 \\
	0&0 \\
	0&1&1
\end{ytableau}
	\hspace{.5cm}
\begin{ytableau}
	\none\\
	0 \\
	0&0 \\
	1&1&1
\end{ytableau}

		\caption{A type-B permutation tableau of size $2$ with $U_2=1$ and the collection of tableaux obtained by extending this tableau to size $3$.}
	\end{center}
\end{figure}

One tableau in $\mathcal{B}_{n-1}$ extends to a group of tableaux in $\mathcal{B}_{n}$ depending on the choices made in the extension--whether a row or a column is inserted and if a column is inserted, whether $0$'s or $1$'s are inserted in the unrestricted rows. Let $\mathcal{F}_{n-1}$ denote the $\sigma$-subalgebra on $\mathcal{B}_{n}$ that groups together all the tableaux obtained by extending the same tableau in $\mathcal{B}_{n-1}$ (see Figure~\ref{extensions}).

Using the terminology and notation described above, several key relationships between $\mathcal{B}_{n}$ and $\mathcal{B}_{n-1}$ were established in \cite{HL} such as the following relationship between the measures,
\begin{equation}\label{eq:exx}
\mathbb{E}_{n}X_{n-1}=\frac{1}{n}\mathbb{E}_{n-1}\left(2^{U_{n-1}}X_{n-1}\right)
\end{equation}
where $X_{n-1}$ is a random variable defined on $\mathcal{B}_{n-1}$. 

The conditional distribution of $U_n$ was also derived, 
\begin{equation}\label{eq:unrestricted}
\mathcal{L}(U_{n}|\mathcal{F}_{n-1})=1+\text{Bin}(U_{n-1})
\end{equation}
where $\text{Bin}(U_{n-1})$ is a binomial random variable with $U_{n-1}$ trials and $1/2$ probability of success.

Recall that  $S_k$ denotes a south step and $W_k$ denotes a west step at the $k$th position. The probability of a south step occurring during the extension procedure was derived in \cite{HL}[Equation 14],
\begin{equation}\label{eq:probsouthstep}
P(S_k|\mathcal{F}_{k-1})=\frac{1}{2^{U_{k-1}+1}}.
\end{equation}

The following calculations can be found in the proofs of \cite{HL}[Proposition~5] and  \cite{HL}[Theorem~7]  and will be established here in order to simplify our work which utilize these results frequently.

\begin{lemma}\label{MeasureBD}
	For any random variable $X$ on $\mathcal{B}_{m-1}$ and $a\in\mathbb{R}$,
	\[
	\mathbb{E}_mXa^{U_m}=\frac{a}{m}\mathbb{E}_{m-1}X\left(a+1\right)^{U_{m-1}}.
	\]
\end{lemma}
\begin{proof}
	Applying Equation~\ref{eq:exx} and the law of total expectation respectively,
		\begin{align*}
	\mathbb{E}_mXa^{U_m}&=\frac{1}{m}	\mathbb{E}_{m-1}X\cdot2^{U_{m-1}}a^{U_m}\\
	&=\frac{1}{m}\mathbb{E}_{m-1}X\cdot2^{U_{m-1}}\mathbb{E}(a^{U_m}|\mathcal{F}_{m-1}).
	\end{align*}
	Now we will apply Equation~\ref{eq:unrestricted} and use the fact that $\mathbb{E} a^{\text{Bin}(n)}=\left(\frac{a+1}{2}\right)^{n}$ to obtain the desired result,
	\begin{align*}
		\mathbb{E}_m(a^{U_m})&=\frac{1}{m}\mathbb{E}_{m-1}X\cdot 2^{U_{m-1}}\mathbb{E}(a^{1+\text{Bin}(U_{m-1})}|\mathcal{F}_{m-1})\\
	&=\frac{1}{m}\mathbb{E}_{m-1}X\cdot 2^{U_{m-1}}a\left(\frac{a+1}{2}\right)^{U_{m-1}}\\
	&=\frac{a}{m}\mathbb{E}_{m-1}X\left(a+1\right)^{U_{m-1}}.
	\end{align*}

\end{proof}

\begin{lemma}\label{DCompose} For any $a\in\mathbb{R}$,
	\[
	\mathbb{E}_m(a^{U_m})=\frac{\Gamma(m+a-1)}{m!\cdot \Gamma(a-1)}.
	\]
\end{lemma}
\begin{proof}
Applying Lemma~\ref{MeasureBD} $m-1$  times,
\[
\mathbb{E}_m(a^{U_m})=\frac{a(a+1)\cdots(a+m-2)}{m\cdots2}\mathbb{E}_1(a+m-1)^{U_1}.
\]
Since $U_1=1$ for all tableau of size $1$,
\[
\mathbb{E}_m(a^{U_m})=\frac{a(a+1)\cdots(a+m-1)}{m\cdots2}=\frac{\Gamma(m+a-1)}{m!\cdot \Gamma(a-1)}.
\]

\end{proof}


\section{Expected number of rows}\label{sectionrows}
In this section, we will determine the expected number of rows in type-B permutation tableaux of size $n$. Note that we are only considering the number of rows in the original Ferrers diagram. This number is equivalent to the number of south steps in the diagram and therefore interesting in terms of the PASEP. Besides, the total number of rows in the shifted Ferrers diagram of a type-B permutation tableaux is simply $n$. This is because the number of rows added in the shifted Ferrer's diagram is equal to the number of columns in the original. 

We will first calculate the probability of a south step at the $k$th position in a type-B permutation tableaux of size $n$ and then use that result to determine the expected number of rows.

\begin{proposition}\label{prop:expsouthstep}

For $1 \le k \le n$, 
\newline
\[P_{n}(I_{S_k})=\frac{1}{2}\left(1-\frac{k-1}{n}\right).\]

\end{proposition}

\begin{proof}
Since we are working with an indicator random variable, 
\[P_{n}(I_{S_k})=E_{n}(I_{S_k}).\]
Applying Lemma~\ref{MeasureBD} $(n-k)$ times,
\begin{align*}
E_{n}(I_{S_k})&=E_{n}(I_{S_k}1^{U_n})\\
&=\frac{(n-k)!}{n\cdot (n-1)\cdots (k+1)}E_{k}(I_{S_{k}}(n-k+1)^{U_{k}}).
\end{align*}
Using the law of total expectation to condition on $\mathcal{F}_{k-1}$,
\begin{align*}
E_{n}(I_{S_k})&=\frac{(n-k)!}{n\cdot (n-1)\cdots (k+1)}E_{k}\mathbb{E}(I_{S_{k}}(n-k+1)^{U_{k}}\,\vert\, \mathcal{F}_{k-1})\\
\end{align*}
Since the kth step is south, $U_{k}=U_{k-1}+1$ and $U_{k-1}$ is constant under the conditional expectation. Thus,
\begin{align*}
E_{n}(I_{S_k})&=\frac{(n-k)!}{n\cdot (n-1)\cdots (k+1)}E_{k}\mathbb{E}(I_{S_{k}}(n-k+1)^{U_{k-1}+1}\,\vert\,\mathcal{F}_{k-1})\\
&=\frac{(n-k+1)!}{n\cdot (n-1)\cdots (k+1)}E_{k}(n-k+1)^{U_{k-1}}\mathbb{E}(I_{S_{k}}\,\vert\,\mathcal{F}_{k-1})\\
&=\frac{(n-k+1)!}{n\cdot (n-1)\cdots (k+1)}E_{k}(n-k+1)^{U_{k-1}}\mathbb{P}(I_{S_{k}}\,\vert\, \mathcal{F}_{k-1})
\end{align*}
where the last equality follows since we are working with an indicator random variable. Applying Equation~\ref{eq:probsouthstep} and Equation~\ref{eq:exx} respectively,
\begin{align}
E_{n}(I_{S_k})&=\frac{(n-k+1)!}{n\cdot (n-1)\cdots (k+1)}E_{k}(n-k+1)^{U_{k-1}}\frac{1}{2^{U_{k-1}+1}}\nonumber\\
&=\frac{(n-k+1)!}{n\cdot (n-1)\cdots (k+1)k}E_{k-1}2^{U_{k-1}}(n-k+1)^{U_{k-1}}\frac{1}{2^{U_{k-1}+1}}\nonumber\\
&=\frac{(n-k+1)!}{n\cdot (n-1)\cdots (k+1)k}\cdot\frac{1}{2}E_{k-1}(n-k+1)^{U_{k-1}}\label{EgOneSStep}.
\end{align}
Applying Lemma~\ref{DCompose} with $m=k-1$ and $a=n-k+1$,
\begin{align*}
E_{n}(I_{S_k})&=\frac{1}{2}\cdot\frac{(n-k-1)!}{n\cdot (n-1)\cdots k}\left(\frac{(n-1)!}{(k-1)!(n-k)!}\right)\\
&=\frac{1}{2}\left(1-\frac{k+1}{n}\right)
\end{align*}
as desired.

\end{proof}

\begin{theorem}\label{thm:exprows}
The expected number of rows in type-B permutation tableaux of size $n$ is $(n+1)/4$. 
\end{theorem}

\begin{proof}
	Let $R_n$ denote the number of rows in type-B permutation tableaux of size $n$. Notice that
	\[\mathbb{E}_{n}R_{n}=\mathbb{E}_{n}\left(\sum_{k=1}^{n}I_{S_k}\right)=\sum_{k=1}^{n}\mathbb{E}_n{}I_{S_k}.\]
	Applying Proposition \ref{prop:expsouthstep},
	\begin{align*}
	\sum_{k=1}^{n}\mathbb{E}_n{}I_{S_k}&=\sum_{k=1}^{n}\frac{1}{2}\left(1-\frac{k-1}{n}\right)\\
	&=\frac{n}{2}-\frac{1}{2n}\left(\frac{n(n+1)}{2}\right)+\frac{1}{2}\\
	&=\frac{n+1}{4}.
	\end{align*}
	
\end{proof}



\section{Expected number of unrestricted rows}\label{sectionunrows}
In this section, we have an immediate result on the expected number of unrestricted rows in type-B permutation tableaux of size $n$. An analogous calculation was performed in \cite{CH} for permutation tableaux and each equation used for that calculation apply in the same manner to type-B permutation tableaux. Therefore, we omit the proof as there is no difference from that of the proof of \cite{CH}[Lemma 4].

\begin{theorem}\label{thm:expunrows}
The expected number of unrestricted rows in type-B permutation tableaux of size $n$ is $H_{n}$, where $H_{n}$ is the $n$th harmonic number.
\end{theorem}


\section{Expected number of ones on the diagonal}\label{sectiondiagonal}
In this section, we will determine the expected number of ones along the diagonal in type-B permutation tableaux of size $n$. For a fixed tableau of size $n$, let $D_{n}$ denote the number of ones on the diagonal, and let $G_k$ denote the position of the topmost one on the $k$th step, numbering only unrestricted rows. For a south step, we leave $G_k$ undefined.

Notice that when $G_k = 1$, the kth diagonal cell contains a $1$ and otherwise, does not. Therefore, 
\[
\mathbb{E}_{n}D_n=\mathbb{E}_{n}\sum_{k=1}^n I_{G_{k}=1}=\sum_{k=1}^n \mathbb{E}_{n}I_{G_{k}=1}.
\] 
To calculate the expectation on the right--hand side, we will use the following property of a binomial random variable  as was used in \cite{CH},
\begin{equation}\label{eq:6}
	\mathbb{E}I_{G=1}a^{G+\text{Bin}(m-G)}=\frac{a}{a+1}\left(\frac{a+1}{2}\right)^{m}.
\end{equation}
In this equation, $G$ represents the position of the first success in the binomial random variable Bin$(m)$.

\begin{proposition}\label{prop:exponesdiagonal}

For $1 \le k \le n$
\[\mathbb{E}_nI_{G_{k}=1}=\frac{1}{2}\]

\end{proposition}

\begin{proof}
	Consider $\mathbb{E}_nI_{G_{k}=1}1^{U_n}$ and apply Lemma~\ref{MeasureBD} $(n-k)$ times,
	\begin{align*}
	\mathbb{E}_nI_{G_{k}=1}&=\frac{(n-k)!}{n\cdot(n-1)\cdots(k+1)}\mathbb{E}_k(I_{G_{k}=1}(n-k+1)^{U_k}).
	\end{align*}
	Now apply the law of total probability to condition on $\mathcal{F}_{k-1}$,
	\begin{align*}
	\mathbb{E}_nI_{G_{k}=1}&=\frac{(n-k)!}{n\cdot(n-1)\cdots(k+1)}\mathbb{E}_k\mathbb{E}(I_{G_{k}=1}(n-k+1)^{U_k}\,|\,\mathcal{F}_{k-1}).
	\end{align*}
	Since we are working with the indicator random variable $I_{G_{k}}=1$, we may assume there is a $1$ on the diagonal at the $k$th step. This implies that the $k$th step is $W$ and therefore, a column is created. For type-B permutation tableaux, there are $U_{k-1}+1$ unrestricted rows in a column created on the $k$th step (see Section~\ref{PApproach} for more details). In our case, one of those rows (the first one) is occupied by a $1$. Therefore, there are $U_{k-1}$ unrestricted rows remaining. Therefore, $U_k=$ Bin$(U_{k-1})+1$ in the conditional expectation,
		\begin{align*}
	\mathbb{E}_nI_{G_{k}=1}&=\frac{(n-k)!}{n\cdot(n-1)\cdots(k+1)}\mathbb{E}_k\mathbb{E}(I_{G_{k}=1}(n-k+1)^{\text{Bin}(U_{k-1})+1}\,|\,\mathcal{F}_{k-1})
	\end{align*}
		Applying Equation~\ref{eq:6} with $m = U_{k-1}+1$, 
		\begin{align*}
		\mathbb{E}_nI_{G_{k}=1}&=\frac{(n-k)!}{n\cdot(n-1)\cdots(k+1)}\mathbb{E}_k\left(\frac{n-k+1}{n-k+2}\left(\frac{n-k+2}{2}\right)^{U_{k-1}+1}\right).
		\end{align*}
		Applying Equation~\ref{eq:exx},
		\begin{align*}
		\mathbb{E}_nI_{G_{k}=1}&=\frac{(n-k+1)(n-k)!}{ 2n\cdot(n-1)\cdots(k+1)k}\mathbb{E}_{k-1}\left(2^{U_{k-1}}\left(\frac{n-k+2}{2}\right)^{U_{k-1}}\right)\\
		&=\frac{(n-k+1)(n-k)!}{ 2n\cdot(n-1)\cdots(k+1)k}\mathbb{E}_{k-1}(n-k+2)^{U_{k-1}}
		\end{align*}
		Finally, apply Lemma~\ref{DCompose} to the expectation with $m = k-1$ and $a = n-k+2$ to obtain the result,
			\begin{align*}
		\mathbb{E}_nI_{G_{k}=1}
		&=\frac{(n-k+1)(n-k)!}{2 n\cdot(n-1)\cdots(k+1)k}\left(\frac{n!}{(k-1)!(n-k+1)!}\right)\\
		&=\frac{1}{2}.
		\end{align*}
		
\end{proof}

\begin{theorem}\label{thm:exponesdiagonal}
	The expected number of ones on the diagonal of a type-B permutation tableau of size $n$ is $n/2$. 
\end{theorem}

\begin{proof} Applying Proposition~\ref{prop:exponesdiagonal},

\[\mathbb{E}_{n}D_n = \sum_{k=1}^{n}\mathbb{E}_{n}I_{G_{k}=1}=\sum_{k=1}^{n}\frac{1}{2}=\frac{n}{2}.\]

\end{proof}

\section{Adjacent south steps}\label{sectiontworows}
In this section, we will determine the expected number of two adjacent south steps in type-B permutation tableaux of size $n$. This statistic corresponds to adjacent filled sites in the PASEP.

First, we will calculate the probability of a south step at both the kth and the (k+1)th position of a type-B permutaiton tableaux of size $n$ and then use that result to compute the expected number of adjacent south steps.  

\begin{proposition}\label{prop:twosouthstep}
	
	For $2 \le k \le n$ ,	
	\[P_{n}(I_{S_k, S_{k-1}})=\displaystyle\frac{(n-k+1)^2}{4n(n-1)}
	\]
	
\end{proposition}
\begin{proof}
	As in Proposition~\ref{prop:expsouthstep}, we will consider the expected value of $I_{S_k, S_{k-1}}$ and work to reduce the measure. The beginning of the proof is almost identical to Proposition~\ref{prop:expsouthstep} until we have reduced to the $(k-1)$th measure, so we will begin from that point (see Equation~\ref{EgOneSStep}).
	\[	
	\mathbb{E}_n(I_{S_k}I_{S_{k-1}})=\frac{(n-k+1)!}{n\cdot (n-1)\cdots (k+1)k}\cdot\frac{1}{2}E_{k-1}I_{S_{k-1}}(n-k+1)^{U_{k-1}}.
	\]
	Using the law of total expectation to condition on $\mathcal{F}_{k-2}$,
	\begin{align*}
	\mathbb{E}_n(I_{S_k}I_{S_{k-1}})&=\frac{(n-k+1)!}{n\cdot (n-1)\cdots (k+1)k}\cdot\frac{1}{2}E_{k-1}\mathbb{E}(I_{S_{k-1}}(n-k+1)^{U_{k-1}}\,\vert\, \mathcal{F}_{k-2})\\
	\end{align*}
	Since the $(k-1)$th step is south, $U_{k-1}=U_{k-2}+1$ and $U_{k-2}$ is constant under the conditional expectation. Thus,
	\begin{align*}
	\mathbb{E}_n(I_{S_k}I_{S_{k-1}})&=\frac{(n-k+1)!}{n\cdot (n-1)\cdots (k+1)k}\cdot\frac{1}{2}E_{k-1}\mathbb{E}(I_{S_{k-1}}(n-k+1)^{U_{k-2}+1}\,\vert\, \mathcal{F}_{k-2})\\
	&=\frac{(n-k+1)\cdot(n-k+1)!}{n\cdot (n-1)\cdots (k+1)k}\cdot\frac{1}{2}E_{k-1}(n-k+1)^{U_{k-2}}\mathbb{E}(I_{S_{k-1}}\,\vert\, \mathcal{F}_{k-2})\\
	&=\frac{(n-k+1)\cdot(n-k+1)!}{n\cdot (n-1)\cdots (k+1)k}\cdot\frac{1}{2}E_{k-1}(n-k+1)^{U_{k-2}}\mathbb{P}(I_{S_{k-1}}\,\vert\, \mathcal{F}_{k-2})
	\end{align*}
	where the last equality follows since we are working with an indicator random variable. Applying Equation~\ref{eq:probsouthstep} and Equation~\ref{eq:exx} respectively,
	\begin{align*}
	\mathbb{E}_n(I_{S_k}I_{S_{k-1}})&=\frac{(n-k+1)\cdot(n-k+1)!}{n\cdot (n-1)\cdots (k+1)k}\cdot\frac{1}{2}E_{k-1}(n-k+1)^{U_{k-2}}\frac{1}{2^{U_{k-2}+1}}\\
	&=\frac{(n-k+1)\cdot(n-k+1)!}{n\cdot (n-1)\cdots k(k-1)}\cdot\frac{1}{2}E_{k-2}2^{U_{k-2}}(n-k+1)^{U_{k-2}}\frac{1}{2^{U_{k-2}+1}}\\
	&=\frac{(n-k+1)\cdot(n-k+1)!}{n\cdot (n-1)\cdots k(k-1)}\cdot\frac{1}{4}E_{k-2}(n-k+1)^{U_{k-2}}.
	\end{align*}
	Applying Lemma~\ref{DCompose} with $m=k-2$ and $a=n-k+1$,
	\begin{align*}
	\mathbb{E}_n(I_{S_k}I_{S_{k-1}})&=\frac{(n-k+1)\cdot(n-k+1)!}{n\cdot (n-1)\cdots k(k-1)}\cdot\frac{1}{4}\left(\frac{(n-2)!}{(k-2)!(n-k)!}\right)\\
	&=\displaystyle\frac{(n-k+1)^2}{4n(n-1)}
	\end{align*}
	as desired.
\end{proof}

\begin{theorem}\label{thm:twosouthstep}
The expected number of two adjacent south steps in type-B permutation tableaux of size $n$ is $\frac{2n-1}{24}$. 
\end{theorem}

\begin{proof}
	Notice that
	
	\[
	\mathbb{E}_n\left(\sum_{k=2}^nI_{S_k,S_{k-1}}\right)=\sum_{k=2}^n\mathbb{E}_nI_{S_k,S_{k-1}}.
	\]
	Applying Proposition~\ref{prop:twosouthstep},
	\begin{align*}
	\mathbb{E}_n\left(\sum_{k=2}^nI_{S_k,S_{k-1}}\right)&=\sum_{k=2}^{n}\frac{(n-k+1)^2}{4n(n-1)}\\
	&=\frac{n^2(n-1)}{4n(n-1)}-\frac{2n}{4n(n-1)}\sum_{k=2}^{n}(k-1)+\frac{1}{4n(n-1)}\sum_{k=2}^{n}(k-1)^2\\
	&=\frac{n}{4}-\frac{1}{2(n-1)}\sum_{j=1}^{n-1}j+\frac{1}{4n(n-1)}\sum_{j=1}^{n-1}j^2\\
	&=\frac{n}{4}-\frac{1}{2(n-1)}\cdot\frac{n(n-1)}{2}+\frac{1}{4n(n-1)}\cdot\frac{n(n-1)(2n-1)}{6}\\
	&=\frac{2n-1}{24}.
	\end{align*}
	
\end{proof}

\section{Adjacent west steps}\label{sectiontwocols}
In this section, we will determine the expected number of two adjacent west steps in type-B permutation tableaux of size $n$. This statistic corresponds to adjacent empty sites in the PASEP.

First, we will calculate the probability of a west step at both the kth and the (k+1)th position of a type-B permutation tableaux of size $n$ and then use that result to compute the expected number of adjacent west steps.

\begin{proposition}\label{2AdjWSteps}
	For $2\leq k \leq n$,
	\[
	\mathbb{P}_n\left(I_{W_k,W_{k-1}}\right)=\frac{k}{n}-\frac{3}{2n}+\frac{(n-k+1)^2}{4n(n-1)}.
	\]	
\end{proposition}
\begin{proof}
Since we are working with an indicator random variable, 
\[P_{n}(I_{W_k,W_{k-1}})=E_{n}(I_{W_k}I_{W_{k-1}}).\]
Applying Lemma~\ref{MeasureBD} $(n-k)$ times,
\begin{align*}
E_{n}(I_{W_k}I_{W_{k-1}})&=E_{n}(I_{W_k}I_{W_{k-1}}1^{U_n})\\
&=\frac{(n-k)!}{n\cdot (n-1)\cdots (k+1)}E_{k}(I_{W_k}I_{W_{k-1}}(n-k+1)^{U_{k}}).
\end{align*}
Using the law of total expectation to condition on $\mathcal{F}_{k-1}$,
\begin{align*}
E_{n}(I_{W_k}I_{W_{k-1}})&=\frac{(n-k)!}{n\cdot (n-1)\cdots (k+1)}E_{k}\mathbb{E}(I_{W_k}I_{W_{k-1}}(n-k+1)^{U_{k}}\,\vert\, \mathcal{F}_{k-1}).\\
\end{align*}
Since $I_{W_{k-1}}$ is constant under the conditional expectation, 
\begin{align}\label{beforeComp}
E_{n}(I_{W_k}I_{W_{k-1}})&=\frac{(n-k)!}{n\cdot (n-1)\cdots (k+1)}E_{k}I_{W_{k-1}}\mathbb{E}(I_{W_k}(n-k+1)^{U_{k}}\,\vert\, \mathcal{F}_{k-1}).
\end{align}
Now use the complement to rewrite 
\begin{align*}
\mathbb{E}(I_{W_k}(n-k+1)^{U_{k}}\,\vert\, \mathcal{F}_{k-1}) &=\mathbb{P}(I_{W_k}(n-k+1)^{U_{k}}\,\vert\, \mathcal{F}_{k-1})\\
&=\mathbb{P}((n-k+1)^{U_{k}}\,\vert\, \mathcal{F}_{k-1})-\mathbb{P}(I_{S_k}(n-k+1)^{U_{k}}\,\vert\, \mathcal{F}_{k-1})\\
&=\mathbb{E}((n-k+1)^{U_{k}}\,\vert\, \mathcal{F}_{k-1})-\mathbb{E}(I_{S_k}(n-k+1)^{U_{k}}\,\vert\, \mathcal{F}_{k-1})
\end{align*}
and combine with Equation~\ref{beforeComp},
\begin{align}\label{firstsplit}
E_{n}(I_{W_k}I_{W_{k-1}})&=\frac{(n-k)!}{n\cdot (n-1)\cdots (k+1)}\Big(E_{k}I_{W_{k-1}}\mathbb{E}((n-k+1)^{U_{k}}\,\vert\, \mathcal{F}_{k-1})\nonumber\\
&\hspace{4cm}-E_{k}I_{W_{k-1}}\mathbb{E}(I_{S_k}(n-k+1)^{U_{k}}\,\vert\, \mathcal{F}_{k-1})\Big).
\end{align}

Now let's consider each part in the parenthesis of Equation~\ref{firstsplit} separately. For the first conditional expectation on the right-hand side, we will apply Equation~\ref{eq:unrestricted} and use the fact that $\mathbb{E} a^{\text{Bin}(n)}=\left(\frac{a+1}{2}\right)^{n}$ to obtain

\begin{align*}
E_{k}I_{W_{k-1}}\mathbb{E}((n-k+1)^{U_{k}}\,\vert\, \mathcal{F}_{k-1})&=E_{k}I_{W_{k-1}}\mathbb{E}((n-k+1)^{1+\text{Bin}(U_{k-1})}\,\vert\, \mathcal{F}_{k-1})\\
&=E_{k}I_{W_{k-1}}(n-k+1)\cdot\frac{(n-k+2)^{U_{k-1}}}{2^{U_{k-1}}}.
\end{align*}
Applying Equation~\ref{eq:exx},
\begin{align*}
E_{k}I_{W_{k-1}}\mathbb{E}((n-k+1)^{U_{k}}\,\vert\, \mathcal{F}_{k-1})&=\frac{n-k+1}{k}E_{k-1}2^{U_{k-1}}I_{W_{k-1}}\cdot\frac{(n-k+2)^{U_{k-1}}}{2^{U_{k-1}}}\\
&=\frac{n-k+1}{k}E_{k-1}I_{W_{k-1}}\cdot(n-k+2)^{U_{k-1}}.
\end{align*}

Next, we apply the law of total probability to condition on $\mathcal{F}_{k-2}$ and use the fact that $I_{S_{k-1}}$ is the complement of $I_{W_{k-1}}$,

\begin{align*}
E_{k}I_{W_{k-1}}\mathbb{E}((n-k+1)^{U_{k}}\,\vert\, \mathcal{F}_{k-1})
&=\frac{n-k+1}{k}E_{k-1}\mathbb{E}\left(I_{W_{k-1}}\cdot(n-k+2)^{U_{k-1}}\,\vert\, \mathcal{F}_{k-2}\right)\\
&=\frac{n-k+1}{k}E_{k-1}\Bigg[\mathbb{E}\left((n-k+2)^{U_{k-1}}\,\vert\, \mathcal{F}_{k-2}\right)\\
&\hspace{3cm}-\mathbb{E}\left(I_{S_{k-1}}(n-k+2)^{U_{k-1}}\,\vert\, \mathcal{F}_{k-2}\right)\Bigg].
\end{align*}

 Next, we will apply  Equation~\ref{eq:unrestricted} and the fact that $\mathbb{E} a^{\text{Bin}(n)}=\left(\frac{a+1}{2}\right)^{n}$ to the first conditional expectation on the right-hand side. For the second conditional expectation, since the $(k-1)$th step is south, $U_{k-1}=U_{k-2}+1$ and then $U_{k-2}$ is constant under the conditional expectation. All together, we obtain
\begin{align*}
E_{k}I_{W_{k-1}}\mathbb{E}((n-k+1)^{U_{k}}\,\vert\, \mathcal{F}_{k-1})
&\\
&\hspace{-3cm}=\frac{n-k+1}{k}E_{k-1}\Bigg[(n-k+2)\frac{(n-k+3)^{U_{k-2}}}{2^{U_{k-2}}}\\
&\hspace{1cm}-(n-k+2)(n-k+2)^{U_{k-2}}\mathbb{E}\left(I_{S_{k-1}}\,\vert\, \mathcal{F}_{k-2}\right)\Bigg].
\end{align*}

By Equation~\ref{eq:probsouthstep},
\[
\mathbb{E}(I_{S_{k-1}}\,\vert\, \mathcal{F}_{k-2})=\mathbb{P}(I_{S_{k-1}}\,\vert\, \mathcal{F}_{k-2})=\frac{1}{2^{U_{k-2}+1}}
\]
and therefore,

\begin{align*}
E_{k}I_{W_{k-1}}\mathbb{E}((n-k+1)^{U_{k}}\,\vert\, \mathcal{F}_{k-1})
&\\
&\hspace{-3cm}=\frac{n-k+1}{k}E_{k-1}\Bigg[(n-k+2)\frac{(n-k+3)^{U_{k-2}}}{2^{U_{k-2}}}\\
&\hspace{1cm}-(n-k+2)(n-k+2)^{U_{k-2}}\frac{1}{2^{U_{k-2}+1}}\Bigg].
\end{align*}

 Applying Equation~\ref{eq:exx} and simplifying,
\begin{align*}
E_{k}I_{W_{k-1}}\mathbb{E}((n-k+1)^{U_{k}}\,\vert\, \mathcal{F}_{k-1})
&\\
&\hspace{-3cm}=\frac{n-k+1}{k(k-1)}E_{k-2}\Bigg[2^{U_{k-2}}(n-k+2)\frac{(n-k+3)^{U_{k-2}}}{2^{U_{k-2}}}\\
&\hspace{1cm}-2^{U_{k-2}}(n-k+2)(n-k+2)^{U_{k-2}}\frac{1}{2^{U_{k-2}+1}}\Bigg]\\
&\hspace{-3cm}=\frac{(n-k+1)(n-k+2)}{k(k-1)}E_{k-2}\Bigg[(n-k+3)^{U_{k-2}}-\frac{1}{2}(n-k+2)^{U_{k-2}}\Bigg].
\end{align*}

The last step for this part is to apply Lemma~\ref{DCompose} to each part of the expected value with $m = k-2$ and $a = n-k+3$ and $m = k-2$ and $a = n-k+2$ respectively,
\begin{align}
E_{k}I_{W_{k-1}}\mathbb{E}((n-k+1)^{U_{k}}\,\vert\, \mathcal{F}_{k-1})
&\nonumber\\
&\hspace{-4cm}=\frac{(n-k+1)(n-k+2)}{k(k-1)}\left[\left(\frac{n!}{(k-2)!(n-k+2)!}\right)-\frac{1}{2}\left(\frac{(n-1)!}{(k-2)!(n-k+1)!}\right)\right]\nonumber\\
&\hspace{-4cm}=\frac{1}{k!}\left[\left(\frac{n!}{(n-k)!}\right)-\frac{1}{2}\left(\frac{(n-1)!(n-k+2)}{(n-k)!}\right)\right]\label{part1}
\end{align}

Now returning to Equation~\ref{firstsplit}, let's compute the remaining term.
	
\[
E_{k}I_{W_{k-1}}\mathbb{E}(I_{S_k}(n-k+1)^{U_{k}}\,\vert\, \mathcal{F}_{k-1}).
\]

Since the $k$th step is south, $U_{k}=U_{k-1}+1$ and then $U_{k-1}$ is constant under the conditional expectation. Therefore,
\[
E_{k}I_{W_{k-1}}\mathbb{E}(I_{S_k}(n-k+1)^{U_{k}}\,\vert\, \mathcal{F}_{k-1})=E_{k}I_{W_{k-1}}(n-k+1)^{U_{k-1}+1}\mathbb{E}(I_{S_k}\,\vert\, \mathcal{F}_{k-1}).
\]
By Equation~\ref{eq:probsouthstep},
\[
E_{k}I_{W_{k-1}}\mathbb{E}(I_{S_k}(n-k+1)^{U_{k}}\,\vert\, \mathcal{F}_{k-1})=\mathbb{E}(I_{S_{k}}\,\vert\, \mathcal{F}_{k-1})=\mathbb{P}(I_{S_{k}}\,\vert\, \mathcal{F}_{k-1})=\frac{1}{2^{U_{k-1}+1}}
\]
and therefore,
\[
E_{k}I_{W_{k-1}}\mathbb{E}(I_{S_k}(n-k+1)^{U_{k}}\,\vert\, \mathcal{F}_{k-1})=E_{k}I_{W_{k-1}}(n-k+1)^{U_{k-1}+1}\frac{1}{2^{U_{k-1}+1}}.
\]
Applying Equation~\ref{eq:exx},
 
 \begin{align*}
E_{k}I_{W_{k-1}}\mathbb{E}(I_{S_k}(n-k+1)^{U_{k}}\,\vert\, \mathcal{F}_{k-1})&=\frac{1}{k} E_{k-1}2^{U_{k-1}}I_{W_{k-1}}(n-k+1)^{U_{k-1}+1}\frac{1}{2^{U_{k-1}+1}}\\
&=\frac{n-k+1}{2k} E_{k-1}I_{W_{k-1}}(n-k+1)^{U_{k-1}}
 \end{align*}
  
Applying the law of total probability and the fact that $I_{S_{k-1}}$ is the complement of $I_{W_{k-1}}$,
 \begin{align*}
E_{k}I_{W_{k-1}}\mathbb{E}(I_{S_k}(n-k+1)^{U_{k}}\,\vert\, \mathcal{F}_{k-2})&=\frac{n-k+1}{2k} E_{k-1}\mathbb{E}\left(I_{W_{k-1}}(n-k+1)^{U_{k-1}}\,\vert\, \mathcal{F}_{k-2}\right)\\
&\hspace{-4cm}=\frac{n-k+1}{2k} E_{k-1}\bigg[\mathbb{E}\left((n-k+1)^{U_{k-1}}\,\vert\, \mathcal{F}_{k-2}\right)-\mathbb{E}\left(I_{S_{k-1}}(n-k+1)^{U_{k-1}}\,\vert\, \mathcal{F}_{k-2}\right)\bigg].
\end{align*}
 Next, we will apply  Equation~\ref{eq:unrestricted} and the fact that $\mathbb{E} a^{\text{Bin}(n)}=\left(\frac{a+1}{2}\right)^{n}$ to the first conditional expectation on the right-hand side. For the second conditional expectation, since the $(k-1)$th step is south, $U_{k-1}=U_{k-2}+1$ and then $U_{k-2}$ is constant under the conditional expectation. All together, we obtain
  \begin{align*}
 E_{k}I_{W_{k-1}}\mathbb{E}(I_{S_k}(n-k+1)^{U_{k}}\,\vert\, \mathcal{F}_{k-1})&\\
 &\hspace{-4cm}=\frac{n-k+1}{2k} E_{k-1}\bigg[(n-k+1)\frac{(n-k+2)^{U_{k-2}}}{2^{U_{k-2}}} -(n-k+1)^{U_{k-2}+1}\mathbb{E}\left(I_{S_{k-1}}\,\vert\, \mathcal{F}_{k-2}\right)\bigg]\\
  &\hspace{-4cm}=\frac{(n-k+1)^2}{2k} E_{k-1}\bigg[\frac{(n-k+2)^{U_{k-2}}}{2^{U_{k-2}}} -(n-k+1)^{U_{k-2}}\mathbb{E}\left(I_{S_{k-1}}\,\vert\, \mathcal{F}_{k-2}\right)\bigg].
 \end{align*}
 
 By Equation~\ref{eq:probsouthstep},

  \begin{align*}
E_{k}I_{W_{k-1}}\mathbb{E}(I_{S_k}(n-k+1)^{U_{k}}\,\vert\, \mathcal{F}_{k-1})&\\
&\hspace{-4cm}=\frac{(n-k+1)^2}{2k} E_{k-1}\bigg[\frac{(n-k+2)^{U_{k-2}}}{2^{U_{k-2}}} -(n-k+1)^{U_{k-2}}\frac{1}{2^{U_{k-2}+1}}\bigg].
\end{align*}
Applying Equation~\ref{eq:exx},
		
	  \begin{align*}
	E_{k}I_{W_{k-1}}\mathbb{E}(I_{S_k}(n-k+1)^{U_{k}}\,\vert\, \mathcal{F}_{k-1})&\\
	&\hspace{-4cm}=\frac{(n-k+1)^2}{2k(k-1)} E_{k-2}\bigg[2^{U_{k-2}}\frac{(n-k+2)^{U_{k-2}}}{2^{U_{k-2}}} -2^{U_{k-2}}(n-k+1)^{U_{k-2}}\frac{1}{2^{U_{k-2}+1}}\bigg]\\
	&\hspace{-4cm}=\frac{(n-k+1)^2}{2k(k-1)} E_{k-2}\left((n-k+2)^{U_{k-2}} -\frac{1}{2}(n-k+1)^{U_{k-2}}\right).
	\end{align*}
	
	The last step for this part is to apply Lemma~\ref{DCompose} to each part of the expected value with $m = k-2$, $a = n-k+2$ and $m = k-2$, $a = n-k+1$ respectively,
	
	 \begin{align}
	E_{k}I_{W_{k-1}}\mathbb{E}(I_{S_k}(n-k+1)^{U_{k}}\,\vert\, \mathcal{F}_{k-1})\nonumber\\
	&\hspace{-4cm}=\frac{(n-k+1)^2}{2k(k-1)} \left[\frac{(n-1)!}{(k-2)!(n-k+1)!} -\frac{1}{2}\frac{(n-2)!}{(k-2)!(n-k)!}\right]\nonumber\\
	&\hspace{-4cm}=\frac{1}{2k!} \left[\frac{(n-1)!(n-k+1)}{(n-k)!} -\frac{1}{2}\frac{(n-2)!(n-k+1)^2}{(n-k)!}\right]\label{part2}.
	\end{align}

	Now returning to Equation~\ref{firstsplit} and plugging in Equation~\ref{part1} and Equation~\ref{part2},

		\begin{align*}
	E_{n}(I_{W_k}I_{W_{k-1}})&=\frac{(n-k)!}{n\cdots (k+1)}\bigg[\frac{1}{k!}\left(\frac{n!}{(n-k)!}-\frac{1}{2}\frac{(n-1)!(n-k+2)}{(n-k)!}\right)\\
	&\hspace{2cm}
	-\frac{1}{2k!}\left(\frac{(n-1)!(n-k+1)}{(n-k)!} -\frac{1}{2}\frac{(n-2)!(n-k+1)^2}{(n-k)!}\right)\bigg]\\
	&=1-\frac{(n-k+2)}{2n}
	-\frac{(n-k+1)}{2n}+\frac{(n-k+1)^2}{4n(n-1)}	\\
	&=\frac{k}{n}-\frac{3}{2n}+\frac{(n-k+1)^2}{4n(n-1)}	
	\end{align*}
	as desired.

\end{proof}

\begin{theorem}\label{thm:twoadjacentw}
	The expected number of two adjacent west steps in type-B permutation tableaux of size $n$ is given by,
	\[\displaystyle\frac{14n-25}{24}+\displaystyle\frac{1}{2n}.\]
\end{theorem}
\begin{proof}
	The expected value can be found by summing the result from Proposition~\ref{2AdjWSteps},
	\begin{align*}
	\sum_{k=2}^{n}\mathbb{P}_n\left(I_{W_k,W_{k+1}}\right)&=\sum_{k=2}^{n} \left(\frac{k}{n}-\frac{3}{2n}+\frac{(n-k+1)^2}{4n(n-1)}\right)\\
	&=\frac{1}{n}\sum_{j=1}^{n-1}(j+1)-\frac{3(n-1)}{2n}+\frac{1}{24n(n-1)}\sum_{j=1}^{n-1}j^2\\
	&=\frac{(n-1)n}{2n}+\frac{n-1}{n}
	-\frac{3(n-1)}{2n}+\frac{(n-1)n(2n-1)}{24n(n-1)}\\
	&=\displaystyle\frac{14n-25}{24}+\displaystyle\frac{1}{2n}.
	\end{align*}

\end{proof}


\section{Conclusion}

\par
In this paper, we calculated the expected value of five statistics on type-B permutation tableaux, specifically rows, unrestricted rows, diagonal ones, adjacent south steps and adjacent west steps. The first three statistics were previously calculated for permutation tableaux in \cite{CH} and the other two calculations are motivated by the PASEP. It would be interesting to compute the expected number of superfluous ones in a type-B permutation tableaux as this was done in \cite{CH} for permutation tableaux but we were unable to do so here for type-B permutation tableaux.



\bibliography{mybib}{}
\bibliographystyle{plain}

\end{document}